\def\ispreprint{}

\documentclass[12pt]{article}

\usepackage{amsmath}
\usepackage{amssymb}
\usepackage{amsthm}

\usepackage[dvips]{graphicx}
\usepackage{color}
\usepackage{epsfig}
\usepackage[mathscr]{eucal}
\newtheorem{theorem}{Theorem}

\newtheorem{prop}[theorem]{Proposition}

\newtheorem{lemma}[theorem]{Lemma}
\newtheorem{cor}[theorem]{Corollary}

\newtheorem{definition}[theorem]{Definition}
\newtheorem{fact}{Fact}

\newcommand{\ft}{\ensuremath{\mathbb{F}_q}}
\newcommand{\ftn}{\ensuremath{\mathbb{F}_{q^n}}}

\newcommand{\floor}[1]{\ensuremath{\left \lfloor {#1} \right \rfloor}}
\newcommand{\ceil}[1]{\ensuremath{\left \lceil {#1} \right \rceil}}

\newcommand{\ZZq}{\mathbb{Z}_q}

\newcommand{\Z}{{\mathbb Z}}

\newcommand{\R}{{\mathbb R}}

\newcommand{\cC}{{\cal C}}
\newcommand{\cE}{{\cal E}}

\newcommand{\cI}{{\cal I}}

\newcommand{\cO}{{\cal O}}

\newcommand{\bP}{{\bf P}}

\newcommand{\KD}{\ensuremath{\widetilde{K}}}
\newcommand{\KL}{\ensuremath{\overleftrightarrow{K}}}

\allowdisplaybreaks

\begin{document}

\ifx\ispreprint\undefined
\title{{\huge \bf Generalized de Bruijn Cycles}}
\author{{\Large Joshua N.\@ Cooper}\\
{\it Department of Mathematics}\\
{\it Courant Institute of Mathematical Sciences, NYU}\\
E-mail: cooper@cims.nyu.edu\\[.25in]
{\Large Ronald L.\@ Graham}\\
{\it Department of Computer Science and Engineering, UCSD}\\
E-mail: graham@ucsd.edu\\}
\else
\title{Generalized de Bruijn Cycles}
\author{Joshua N. Cooper \\ \small{Department of Mathematics} \\ \small{Courant Institute of Mathematical Sciences, NYU} \\[.25in] Ronald L. Graham \\ \small{Department of Computer Science and Engineering, UCSD}}
\fi

\date{\today}
\maketitle
\ifx\ispreprint\undefined 

\pagebreak \fi

\begin{abstract}
For a set of integers $\cI$, we define a $q$-ary $\cI$-cycle to be a assignment of the symbols $1$ through $q$ to the integers modulo $q^n$ so that every word appears on some translate of $\cI$.  This definition generalizes that of de Bruijn cycles, and opens up a multitude of questions.  We address the existence of such cycles, discuss ``reduced'' cycles (ones in which the all-zeroes string need not appear), and provide general bounds on the shortest sequence which contains all words on some translate of $\cI$.  We also prove a variant on recent results concerning decompositions of complete graphs into cycles and employ it to resolve the case of $|\cI|=2$ completely.
\end{abstract}

\section{Introduction} \label{section:introduction}

A \textit{de Bruijn cycle} of order $n$ is a $q$-ary sequence $(S(0),\ldots,S(q^n-1))$ so that every $q$-ary $n$-word appears in a ``window'' $(S(j),\ldots,S(j+n-1))$ for some $j$ (indices taken modulo $q^n$).  A \textit{reduced} de Bruijn cycle is a string of length $q^n-1$ which achieves every $n$-word in some window, except for the word $0^n$.  In this paper, we are concerned with such objects when the notion of ``window'' is generalized.

Given a sequence $\mathcal{I} = \{i_j\}_{j=1}^n \subset \mathbb{Z}_{q^n}$, we say that the map $\chi: \Z_{q^n} \rightarrow [q]$ (resp., $\chi: \Z_{q^n-1} \rightarrow [q]$) is an $\mathcal{I}$-cycle ($\mathcal{I}^\ast$-cycle) if, for every word $W \in [q]^n$ (resp., every word $W \in [q]^n \setminus \{0^n\}$), there exists a $t \in \mathbb{Z}_{q^n}$ (resp., $t \in \mathbb{Z}_{q^n-1}$) so that $\chi(i_j+t) = W(j)$.  If such a sequence exists for $\cI$, we say that $\cI$ is $q$-valid (resp., $q^\ast$-valid).  We will often abuse notation by writing $\chi(\cI)$ for the $n$-word $(\chi(i_1),\ldots,\chi(i_n))$.  Furthermore, we {\it also} use ``cycle'' to refer to sequences of edges in a directed graph each of whose tail is the head of the previous one, and which return to their starting point.  It should always be clear from context which of these definitions we intend -- though, often, the notions will coincide!

It is classical that an $\cI$-cycle and an $\cI^\ast$-cycle exist for $\cI = \{1,\ldots,n\}$.  In Section \ref{sec:unreduced}, we address the validity of other sets $\cI$.  As it turns out, the question is rather difficult in general, and we solve the problem completely only for sets of cardinality $2$.  We present some general constructions and a number of computational results, and we discuss a graph-theoretic question whose solution is equivalent to the case of $\cI$ being an arithmetic progression.  In the next section, we discuss the existence of reduced de Bruijn cycles, and have a greater degree of success in characterizing the $q^\ast$-valid sets.  Then, in Section \ref{sec:approx}, the issue of ``approximate'' cycles is discussed, and we present a nearly optimal bound on their length.  The following section contains a proof of a graph-theoretic decomposition result that is used in Section \ref{sec:unreduced} and which solves a variant of a family of problems that has appeared recently in the literature.  We finish with a number of open questions and suggestions for future investigation.

\section{Unreduced Cycles} \label{sec:unreduced}

It is easy to see that, to determine the two-element $q$-valid sets, we need only examine the sets $\{0,d\}$ with $d | q^2$.  Indeed, for any $k \in \Z_q^\times$, if there exists an $\cI$-cycle $\chi$ for $\{0,d\}$, then $\chi^\prime : s \rightarrow \chi(k^{-1}s)$ is an $\cI$-cycle for $\{0,kd\}$.  In addition, it is clear that the validity of $\cI$ is equivalent to the validity of $\cI+b$ for any $b \in \Z_q$.  The same arguments apply to sets $\cI$ whose elements are in arithmetic progression: we need only examine the cases when the difference $d$ divides $q^n$.

Suppose, then, that $d | q^n$, and consider $D_q^n$, the $n^\textrm{th}$ $q$-ary {\it de Bruijn digraph}, i.e., the digraph whose vertices are the $q$-ary $n$-strings and which has an edge from $x$ to $y$ if the last $n-1$ symbols of $x$ are the same as the first $n-1$ symbols of $y$.  (Note that some edges have loops attached.)  Then the set $AP(n,d) = \{0,d,2d,\ldots,(n-1)d\}$ is $q$-valid iff there is a partition of the edges of $D_q^{n-1}$ into $d$ cycles each of length $q^n/d$, because we may write $\chi(j)=C_{a}(b)$, where $C_{a}$ is the $a^\textrm{th}$ such cycle and $j = ad+b$.

Using $D_q^1$, which is simply a complete directed graph with loops on $q$ vertices, we may state a condition equivalent to the $q$-validity of $\{0,d\}$: that $D_q^1$ is fully decomposable into cycles of length $q^2/d$.  Theorem \ref{trails}, which appears Section \ref{sec:trails}, says that this is possible precisely when $q^2/d > 2$.  Therefore, we have

\begin{theorem} There exists an $\cI$-cycle for $\cI=\{0,d\}$ if and only if $q^2/d \neq 2$.
\end{theorem}

The situation for sets with $|\cI| > 2$ appears significantly more complicated, even for arithmetic progressions.  However, the above invalidity result for $d = q^2/2$ has an immediate analogue for larger $n$:

\begin{prop} For any $r$, $q$, with $r | q$, the set $AP(r,q^r/r)$ is $q$-invalid.
\end{prop}
\begin{proof} Suppose an $\cI$-cycle $\chi$ existed.  Then, we may assume without loss of generality that $(\chi(0),\ldots,\chi((r-1)q^n/r))$ is the all-zeroes vector.  But, then $(\chi(q^r/r),\chi(2q^r/r),\ldots,\chi(0))$ is also, a contradiction.
\end{proof}

On the other hand, we can construct a large family of $AP(n,q)$-cycles.  Form the \emph{quotient graph} $G(n)$ from $D_{q}^{n}$ by identifying two vertices $x$ and $x^\prime$ of $D_{q}^{n}$ if $ x - x^\prime = k$ for some $k \in \ZZq$, i.e., $x_{i} - x'_{i} = k$ for $1 \leq i \leq n$.

\begin{fact}
$G(n)$ is isomorphic to $D_{q}^{n-1}$.
\end{fact}

\proof
For $x = (x_{1}, x_{2}, \ldots ,x_{n})$, consider the map
$$
\lambda \colon x \mapsto (x_{1}-x_{2}, x_{2}-x_{3}, \ldots ,x_{n-1}-x_{n}).
$$
It is easy to check that $\lambda$ is well-defined on $G(n)$, invertible and preserves directed edges (i.e., $(x,y)$ is an edge in $D_{q}^{n}$ if and only if $(\lambda(x), \lambda(y))$ is an edge in $G(n)$). Note that the inverse map $\lambda^{-1}$ doesn't necessarily preserve \emph{cycles}, though. However, it is not hard to show the following. Suppose that $(q_{0}, q_{1}, \ldots ,q_{r-1})$ is a cycle in $G(n) \cong D_{q}^{n-1}$, i.e.,
$$
((q_{i+1},q_{i+2}, \ldots ,q_{i+n-1}), (q_{i+2},q_{i+3}, \ldots , q_{i+n}))
$$
is an edge for all $i$, where the indices are reduced modulo $r$. Then this cycle ``lifts'' under $\lambda^{-1}$ to a cycle in $D_{q}^{n}$ if and only if $\sum_{i=0}^{r-1} q_{i} = 0 \pmod{q}$. 

Observe that any Eulerian cycle $C$ in $D_{q}^{n-1}$ satisfies this sum condition.  Hence, it lifts to a cycle $C^{+}$ in $D_{q}^{n}$ going through exactly \emph{one} of the $q$ points in each equivalence class.  In fact, we can form $q$ disjoint cycles $\{C_{j}^{+}\}_{j=1}^q$ from this cycle $C^{+}$ by translating each point in it by some fixed constant $q \in \ZZq$. 

Finally, we can form a cyclic sequence $S$ in $D_{q}^{n}$ containing all of its vertices by ``splicing together'' these $q$ cycles $C_{j}^{+}$ in the obvious way. Since $C$ was in fact a de Bruijn cycle for $(n-1)$-tuples, then it easily checked that $S$ is an $\mathcal{I}$-cycle with $\mathcal{I} = \{0, q, 2q, \ldots ,(n-1)q\}$.

As an example, consider the cycle $001122021$ for $n=2, q=3$.  We can lift this to $100021200$ for $n=3, q=3$, form the two translates $211101011$ and $022210122$ and splice them together to get $021210210210102021102210210$, which is a $\{0,3,6\}$-cycle.

In fact, since there are many ways of choosing the first de Bruijn cycle ($[(q-1)!]^{q^{n-1}} \cdot q^{q^{n-2}-n+1}$, to be precise), and many ways of splicing them together ($q!$), there quite a few ways of producing such cycles.  Unfortunately, it is not possible to iterate this construction, since the cycles $C_{q}^{+}$ do not themselves have the zero-sum property needed to lift them again.

\begin{prop} $AP(n,q)$ is $q$-valid for any $n$, $q$.
\end{prop}

Suppose that $\{0,d,2d\}$ is $q$-valid, where $8 | k = q^3 / d$.  Then there exists a decomposition of the edges of $D^2_q$ into cycles of length $k$.  Write $\cE$ for the set of even symbols in $[2q]$ and $\cO$ for the set of odd symbols.   We may think of $D_{2q}^n$ as being composed of four parts: $U_1 = \cE \times \cE$, $U_2 = \cE \times \cO$, $U_3 = \cO \times \cO$, and $U_4 = \cO \times \cE$.  $U_2$ and $U_4$ contain no edges; $U_1$ and $U_3$ are copies of $D^2_q$.  We may therefore decompose $U_1$ and $U_3$ into $k$-cycles.  The remaining edges may be decomposed into $4$-cycles and $8$-cycles as follows.  For each pair $(a,b),(b,c)$ with $a,b \in \cE$ and $c \in \cO$, define a cycle $\cC(a,b,c) = ((a,b),(b,c),(c,a+b+c),(a+b+c,a))$, with addition modulo $2q$.  The resulting $4$-cycles partition all edges which do not belong to $U_1 \times U_1$, $U_3 \times U_3$, $U_2 \times U_4$, or $U_4$ times $U_2$.  The edges in these final two classes come in pairs $\{(b,c),(c,b)\}$.  We may attach $\{(b,c),(c,b)\}$ and $\{(2b+c+2,b),(b,2b+c+2)\}$ to the cycle $\cC(b+2,b,c)$ for each $b$ even and $c = 1 \pmod{4}$, thus turning it into a $8$-cycle.  Doing so accounts for all the remaining edges exactly once.

The result is a set of $8$-cycles and $4$-cycles.  We may partition them into classes so that each class has exactly $k$ edges, and ``join'' each class at $U_1$ into a cycle of length $k$.  The resulting decomposition of $D_{2q}^n$ gives rise to a $\{0,8d,16d\}$-cycle.  Hence, we have the following.

\begin{prop} If $AP(3,d)$ is $q$-valid, where $8 | q^3/d$, then $AP(3,8d)$ is $2q$-valid.
\end{prop}

\begin{cor} $AP(3,8^k)$ is $r2^{k+1}$-valid for all $k \geq 0$ and $r \geq 1$.
\end{cor}
\begin{proof} Begin with the $2r$-valid $\{0,1,2\}$ and iterate the above proposition.
\end{proof}

To illustrate the complicated nature of the $d > 2$ case, we offer the following computational observations.  By ``affine equivalence'', we mean a map $\sigma : s \rightarrow ks + b$ for some $k \in \Z_{q^n}^\times$ and $b \in \Z_{q^n}$.  Clearly, the partition of index sets into valid and invalid is refined by the partition into affine equivalence classes.

\begin{enumerate}
\item For $(q,n)=(3,3)$, the only invalid index sets are $\{k,k+9,k+18\}$ for $k=0\ldots8$.
\item The only $2$-invalid $3$-set (up to affine equivalence) is $\{0,1,3\}$.  For $(q,n) = (2,4)$, the only valid $\mathcal{I}$'s (up to affine equivalence) are the nine sets $\{0,1,2,3\}$, $\{0,1,2,6\}$, $\{0,1,2,7\}$, $\{0,1,3,4\}$, $\{0,1,3,7\}$, $\{0,1,3,8\}$, $\{0,1,3,9\}$, $\{0,1,3,14\}$, and, of course, $\{0,2,4,6\}$.
\item For $(q,n)=(2,5)$, the following list contains one representative of each equivalence class of invalid index sets, in lexicographic order:
\end{enumerate}

\begin{center}
\begin{tabular}{|c|c|c|c|c|}
\hline
0,1,2,3,12 & 0,1,2,4,7 & 0,1,2,4,9 & 0,1,2,4,12 & 0,1,2,4,23 \\
0,1,2,4,24 & 0,1,2,4,25 & 0,1,2,4,26 & 0,1,2,4,27 & 0,1,2,5,7 \\
0,1,2,5,8 & 0,1,2,5,9 & 0,1,2,5,13 & 0,1,2,5,14 & 0,1,2,5,15 \\
0,1,2,5,16 & 0,1,2,5,19 & 0,1,2,5,20 & 0,1,2,5,21 & 0,1,2,5,22 \\
0,1,2,5,24 & 0,1,2,5,25 & 0,1,2,5,26 & 0,1,2,6,9 & 0,1,2,6,11 \\
0,1,2,6,13 & 0,1,2,6,14 & 0,1,2,6,15 & 0,1,2,6,16 & 0,1,2,6,17 \\
0,1,2,6,19 & 0,1,2,6,21 & 0,1,2,6,23 & 0,1,2,6,25 & 0,1,2,6,26 \\
0,1,2,7,11 & 0,1,2,7,14 & 0,1,2,7,15 & 0,1,2,7,19 & 0,1,2,7,22 \\
0,1,2,7,23 & 0,1,2,7,24 & 0,1,2,8,12 & 0,1,2,8,13 & 0,1,2,8,16 \\
0,1,2,8,17 & 0,1,2,8,18 & 0,1,2,8,19 & 0,1,2,8,23 & 0,1,2,8,24 \\
0,1,2,8,25 & 0,1,2,9,13 & 0,1,2,9,14 & 0,1,2,9,15 & 0,1,2,9,16 \\
0,1,2,9,17 & 0,1,2,9,19 & 0,1,2,9,20 & 0,1,2,9,21 & 0,1,2,9,22 \\
0,1,2,9,25 & 0,1,2,10,14 & 0,1,2,10,15 & 0,1,2,10,16 & 0,1,2,10,17 \\
0,1,2,10,18 & 0,1,2,10,20 & 0,1,2,11,13 & 0,1,2,11,14 & 0,1,2,11,15 \\
0,1,2,11,16 & 0,1,2,11,19 & 0,1,2,12,15 & 0,1,2,12,17 & 0,1,2,12,19 \\
0,1,2,13,16 & 0,1,2,13,18 & 0,1,2,13,20 & 0,1,2,14,17 & 0,1,2,15,18 \\
\hline
\end{tabular}

\begin{tabular}{|c|c|c|c|c|}
\hline
0,1,2,16,18 & 0,1,3,4,9 & 0,1,3,4,11 & 0,1,3,4,12 & 0,1,3,4,15 \\
0,1,3,4,16 & 0,1,3,5,9 & 0,1,3,5,11 & 0,1,3,5,12 & 0,1,3,5,13 \\
0,1,3,5,15 & 0,1,3,5,17 & 0,1,3,5,21 & 0,1,3,5,22 & 0,1,3,5,24 \\
0,1,3,5,25 & 0,1,3,5,26 & 0,1,3,7,9 & 0,1,3,7,11 & 0,1,3,7,12 \\
0,1,3,7,15 & 0,1,3,7,16 & 0,1,3,7,17 & 0,1,3,7,19 & 0,1,3,7,23 \\
0,1,3,7,24 & 0,1,3,7,27 & 0,1,3,7,30 & 0,1,3,8,10 & 0,1,3,8,12 \\
0,1,3,8,14 & 0,1,3,8,16 & 0,1,3,8,17 & 0,1,3,8,19 & 0,1,3,8,20 \\
0,1,3,8,21 & 0,1,3,8,23 & 0,1,3,8,24 & 0,1,3,9,13 & 0,1,3,9,16 \\
0,1,3,9,17 & 0,1,3,9,20 & 0,1,3,9,25 & 0,1,3,9,26 & 0,1,3,9,28 \\
0,1,3,9,30 & 0,1,3,10,12 & 0,1,3,10,13 & 0,1,3,10,14 & 0,1,3,10,15 \\
0,1,3,10,16 & 0,1,3,10,20 & 0,1,3,10,23 & 0,1,3,10,30 & 0,1,3,12,13 \\
0,1,3,12,16 & 0,1,3,12,24 & 0,1,3,12,25 & 0,1,3,12,27 & 0,1,3,12,28 \\
0,1,3,13,15 & 0,1,3,13,21 & 0,1,3,13,22 & 0,1,3,13,25 & 0,1,3,13,27 \\
0,1,3,13,28 & 0,1,3,14,15 & 0,1,3,15,16 & 0,1,3,15,20 & 0,1,3,15,21 \\
0,1,3,15,22 & 0,1,3,15,23 & 0,1,3,15,25 & 0,1,3,15,27 & 0,1,3,15,28 \\
0,1,3,16,17 & 0,1,3,16,19 & 0,1,3,16,21 & 0,1,3,16,23 & 0,1,3,16,25 \\
0,1,3,16,27 & 0,1,3,17,19 & 0,1,3,17,21 & 0,1,3,17,23 & 0,1,3,17,25 \\
0,1,3,17,27 & 0,1,3,17,28 & 0,1,3,19,27 & 0,1,3,20,24 & 0,1,3,21,24 \\
0,1,3,21,25 & 0,1,3,22,24 & 0,1,3,22,25 & 0,1,3,23,24 & 0,1,3,23,25 \\
0,1,3,23,27 & 0,1,3,24,28 & 0,1,3,25,27 & 0,1,3,27,28 & 0,1,4,5,13 \\
0,1,4,6,12 & 0,1,4,6,14 & 0,1,4,6,17 & 0,1,4,6,18 & 0,1,4,6,20 \\
0,1,4,8,23 & 0,1,4,9,15 & 0,1,4,9,16 & 0,1,4,9,17 & 0,1,4,9,20 \\
0,1,4,12,14 & 0,1,4,12,18 & 0,1,4,13,14 & 0,1,4,14,17 & 0,1,4,14,28 \\
0,1,4,14,29 & 0,1,4,15,16 & 0,1,4,15,20 & 0,1,4,15,23 & 0,1,4,15,28 \\
0,1,4,16,26 & 0,1,4,17,26 & 0,1,4,18,26 & 0,1,4,26,28 & 0,1,5,7,16 \\
0,1,5,8,16 & 0,1,6,8,17 & 0,1,7,8,17 & 0,1,7,9,15 & 0,1,7,9,16 \\
0,1,7,9,17 & 0,1,7,15,16 & 0,1,8,16,17 & 0,1,8,16,24 & 0,2,4,8,14 \\
0,2,4,10,14 & 0,2,4,10,18 & 0,2,4,10,24 & 0,2,4,10,26 & 0,2,4,16,20 \\
0,2,6,18,22 & 0,2,8,16,18 & 0,2,8,16,24 & 0,4,8,16,24 & \\
\hline
\end{tabular}
\end{center}

\section{Reduced Cycles}

Although the definition of $q^\ast$-validity certainly makes sense when $q$ is not a prime power, we restrict our attention to that case in this section.  Therefore, consider $q \geq 2$ a fixed prime power, and take our alphabet to be $\ft$. Let $\alpha$ be a generator of the multiplicative group of the finite field $\ftn$.  Denote by $\mathcal{E}$ the elementary basis for $\ft^n$ over $\ft$.  Given a basis $\mathcal{B} = \{b_1,\ldots,b_n\}$ of $\ftn$ over $\ft$ and an element $\gamma \in \ftn$, write $f_{\mathcal{B}}(\gamma)$ for the element of $\ft^n$ whose $j^\textrm{th}$ coordinate is the coefficient of $b_j$ in the $\mathcal{B}$-representation of $\gamma$.  Then, given a nonzero vector $\mathbf{v} \in \ft^n$, define $\Lambda(\alpha,\mathcal{B},\mathbf{v})$ to be the string whose $j^\textrm{th}$ coordinate (i.e., $\Lambda_j(\alpha,\mathcal{B},\mathbf{v})$, $0 \leq j \leq q^n-2$) is $\mathbf{v}^\intercal f_{\mathcal{B}}(\alpha^j)$.

It is well known that, when $\mathcal{B}=\{\alpha^j : 0 \leq j \leq n-1\}$ and $\mathbf{v}$ has only one nonzero coordinate, $\Lambda(\alpha,\mathcal{B},\mathbf{v})$ is a reduced de Bruijn cycle of order $n$ (e.g., \cite{Fred82}.)  We generalize this result as follows.

\begin{prop} \label{hom} Let $\mathcal{I} = \{i_j\}_{j=1}^n$ be a sequence of distinct integers.  Fix a basis $\mathcal{B}$ of $\ft^n$ over $\ft$, a generator $\alpha \in \ftn^\times$, and a vector $\mathbf{v} \in \ft^n$, and write $\Phi(t)$ for the vector
$$
(\Lambda_{i_1+t}(\alpha,\mathcal{B},\mathbf{v}),\ldots,\Lambda_{i_n+t}(\alpha,\mathcal{B},\mathbf{v}))^\intercal \in \ft^n
$$
with indices taken modulo $q^n-2$.  If the minimal polynomial of $\alpha$ is not a divisor of $\sum_{j=1}^{n} c_j x^{i_j}$ for any nonzero $(c_1,\ldots,c_n)$, then the map $\Psi$ which sends $0$ to $0$ and $\alpha^t$ to $\Phi(t)$ is an isomorphism from the additive group of $\ftn$ to $\ft^n$.
\end{prop}
\begin{proof} First, we show that $\Psi$ is linear.  Write $e_j$ for the elementary $n$-vector whose coordinates are all zero except for a $1$ in the $j^\textrm{th}$ coordinate.  We denote by $M_{\gamma,\mathcal{B}}$ the matrix representing multiplication by $\gamma \in \ftn$ in the $\mathcal{B}$ basis.  It is easy to see that
$$
\Lambda_{k}(\alpha,\mathcal{B},\mathbf{v}) = \mathbf{v}^\intercal f_{\mathcal{B}}(\alpha^{k})
$$
and therefore that
\begin{equation} \label{matrixformula}
\Psi(\gamma) = \sum_{j=1}^{n} e_{j} \mathbf{v}^\intercal f_{\mathcal{B}}(\alpha^{i_j} \gamma) = \sum_{j=1}^{n} e_{j} \mathbf{v}^\intercal M_{\alpha,\mathcal{B}}^{i_j} f_{\mathcal{B}}(\gamma),
\end{equation}
which is obviously linear.

Now, suppose that $\Psi(\gamma)=0$ and $\gamma = \alpha^t$.  If we denote by $S$ the subspace of $\ft^n$ orthogonal to $\mathbf{v}$, then we have $\alpha^{i_j+t} \in f_{\mathcal{B}}^{-1}(S)$ for each $j$.  However, $f_{\mathcal{B}}$ is linear and has a trivial kernel, so all the $\alpha^{i_j+t}$ lie in a subspace of $\ft^n$ of dimension $n-1$ and are therefore linearly dependent.  Since $M_{\alpha,\mathcal{B}}$ is nonsingular, this implies that $\{\alpha^{i_j}\}_{j=1}^{n}$ is a dependent set.  But then we have
$$
\sum_{j=1}^{n} c_j \alpha^{i_j} = 0
$$
for some nonzero $(c_1,\ldots,c_n)$, a contradiction.
\end{proof}

The map $\gamma \mapsto M_{\gamma,\mathcal{B}}$ is actually an isomorphism of fields.  The image is a set of matrices which form a field, i.e., a \textit{matrix field}.  These objects have been studied extensively and thoroughly characterized when the matrices take their entries from a finite field (\cite{Beard86}).

\begin{cor} If the minimal polynomial of $\alpha$, a multiplicative generator of $\mathbb{F}^\times_{q^n}$, is not a divisor of $\sum_{j=1}^{n} c_j x^{i_j}$ for any nonzero $(c_1,\ldots,c_n)$, then $\Lambda(\alpha,\mathcal{B},\mathbf{v})$ is an $\mathcal{I}^\ast$-cycle.
\end{cor}
\begin{proof} By the above argument, $\Lambda(\alpha,\mathcal{B},\mathbf{v})$ contains all nonzero $n$-strings in shifted copies of the index set $\mathcal{I}$.
\end{proof}

We require another definition.

\begin{definition} The index set $\mathcal{I} = \{i_j\}_{j=1}^n$ is called exceptional for $q$ if, for every primitive polynomial $g \in \ft[x]$ of degree $n$, there exists a nonzero vector $(c_1,\ldots,c_n) \in \ft^n$ so that $g$ divides
\begin{equation} \label{rotation}
\sum_{j=1}^{n} c_j x^{i_j}.
\end{equation}
Equivalently, if for every $\alpha$ a generator of $\ftn^\times$, the set $\{\alpha^{i_j}\}$ is linearly dependent over $\ft$, then $\mathcal{I}$ is exceptional for $q$.  An index set which is not exceptional is called ordinary for $q$.
\end{definition}

Note that the exponents in (\ref{rotation}) can be thought of as belonging to $\mathbb{Z}_{q^n-1}$.

\begin{prop} $\mathcal{I}$ is $q^\ast$-valid whenever $\mathcal{I}$ is ordinary for $q$.
\end{prop}
\begin{proof} If $\mathcal{I} = \{i_j\}_{j=1}^n$ is not exceptional, then there exists a primitive polynomial $g \in \ft[x]$ of degree $n$ so that, for all nonzero $(c_1,\ldots,c_n)$, $g$ is not a divisor of $\sum_{j=1}^{n} c_j x^{i_j}$.  Since $x$ is not a root of $\sum_{j=1}^{n} c_j x^{i_j}$ in $\ft[x]/g$, but it is a multiplicative generator of this field, $\Lambda(x,\mathcal{B},\mathbf{v})$ is an $\mathcal{I}^\ast$-cycle for any $\mathcal{B}$ and $\mathbf{v}$.
\end{proof}

Which index sets are ordinary?  We argue that if $(q^n-1,d)=1$, then $(a,a+d,a+2d,\ldots,a+(n-1)d)$ is ordinary.  It is clear that, if $\mathcal{I}$ is ordinary, all of its translates are as well.  We may therefore assume that $a=0$.  Then some irreducible polynomial $g$ of degree $n$ divides $\sum_{j=0}^{n-1} c_j x^{jd}$ for each $(c_1,\ldots,c_n) \neq 0$.  Let $\alpha$ be a root of $g$, so $\alpha^d$ is a root of $\sum_{j=0}^{n-1} c_j x^{j}$.  That this polynomial has degree less than $n$ contradicts the fact that $\alpha^d$ is a generator of $\ftn^\times$.

It is trivial that $\mathcal{I}$ is ordinary if it is a singleton.  If $\mathcal{I}$ has two elements, then it is easy to see that $\mathcal{I}=\{i,j\}$ is ordinary if $q+1 \nmid i-j$, since then $\alpha^{i-j} \not \in \mathbb{F}_q$.  Indeed, the copy of $\mathbb{F}^\times_q$ lying inside of $\mathbb{F}^\times_{q^2}$ is the set $\{\alpha^{k(q+1)}\}_{k=1}^{q-1}$ for any generator $\alpha$.  Conversely, if $q+1 | i-j$, then, for every generator $\alpha$, we have $\alpha^i = c \alpha^j$ for some $c \in \mathbb{F}_q$.  Therefore, a two-element set is ordinary for $q$ if and only if the difference of the elements is not a multiple of $q+1$.

For a prime $p$ and a positive integer $n$, define the \textit{Jacobi logarithm} as follows: for a generator $\alpha$ of $\mathbb{F}_{p^n}^\times$, define $L_\alpha : \mathbb{Z}_{p^n-1} \setminus \{s\} \mapsto \mathbb{Z}_{p^n-1} \setminus \{0\}$ by $1 + \alpha^t = \alpha^{L_\alpha(t)}$, where $s = (p^n-1)/2$ if $p > 2$ and $s=0$ otherwise.

\begin{prop} A three-element set $\{i,j,k\}$ is exceptional if and only if either:
\begin{enumerate}
\item $Q | j-i$,
\item $Q | k-j$,
\item $Q | k-i$,
\item \label{four} or, for all $m$ with $(m,q^3-1)=1$, there exists an $a$ so that
$$
L_\alpha(aQ+m(j-i)) = m(k-i) \!\!\!\mod Q,
$$
\end{enumerate}
where $Q = q^2 + q + 1$, and $\alpha$ is any (fixed) generator of $\mathbb{F}^\times_{q^3}$.
\end{prop}
\begin{proof} There are only two ways that $\alpha^i$, $\alpha^j$, and $\alpha^k$ can be linearly dependent.  Either one of them is an $\mathbb{F}_q$-multiple of another, or, for some triple $\{c_1,c_2,c_3\}$, with $c_i \neq 0$ for all $i$,
\begin{equation}
\label{lindep} c_1 \alpha^i + c_2 \alpha^j + c_3 \alpha^k = 0.
\end{equation}
The former case is precisely the divisibility conditions stated above.  To see that the latter situation is equivalent to condition \ref{four}, we may rewrite (\ref{lindep}) without loss of generality as
\begin{equation}
\label{lindep2} c_4 \alpha^{i-k} + c_5 \alpha^{j-k} = 1
\end{equation}
with $c_4,c_5 \in \mathbb{F}_q$.  Suppose $\{i,j,k\}$ falls into this case, i.e., (\ref{lindep2}) has a solution in $c_4$ and $c_5$.  Since we have assumed that neither $c_4$ nor $c_5$ is zero, we may express each of them in terms of $\alpha$: respectively, $\alpha^{sQ}$ and $\alpha^{tQ}$, for some $s,t$ integers.  Then we may rewrite (\ref{lindep2}) as
$$
\alpha^{sQ+i-k} (1 + \alpha^{Q(t-s)+j-i} ) = \alpha^{sQ+i-k+L_\alpha(Q(t-s)+j-i)} = \alpha^0,
$$
which is to say, that $sQ+i-k+L_\alpha(Q(t-s)+j-i) = 0 \!\!\mod q^3-1$.  (Note that the fact that $L_\alpha$ is not defined on all of $\mathbb{Z}_{q^3-1}$ is not problematic, since the left hand side of (\ref{lindep2}) cannot be zero.)  This equation has a solution in $s$ and $t$ iff there exists an $a$ so that $L_\alpha(aQ+j-i) = k-i \!\!\mod Q$.

If we choose any other generator $\beta$ of $\mathbb{F}_{q^3}^\times$, there is some $m$ such that $(m,q^3-1)=1$ and $\alpha^m = \beta$.  Then $L_\beta(a)=b$ iff $L_\alpha(am)=bm$, so $\{i,j,k\}$ is exceptional iff $Q | j-i$, $Q | k-j$, $Q | k-i$, or, for some $m$ such that $(m,q^3-1)=1$, there exists an $a$ so that
$$
L_\alpha(aQ+m(j-i)) = m(k-i) \!\!\!\mod Q.
$$
\end{proof}

\section{Approximate Cycles} \label{sec:approx}

Since the question of whether an $\cI$-cycle exists appears difficult in general, we may ask instead whether it is possible to find an ``approximate'' $\cI$-cycle.  This question comes in two forms for an index set $\cI = \{i_j\}_{j=1}^n \subset \Z$:
\begin{enumerate}
\item What is the least $N$ for which a $\chi : \Z_N \rightarrow [q]$ exists so that, for every word $W$, there exists an $m$ with $W = \chi(\cI+m)$?
\item What is the least $N$ for which there exists a $\chi$ so that all but $o(q^N)$ words appear as $\chi(\cI+m)$ for some $m$?
\end{enumerate}
Call the former object a {\it Type I approximate cycle} and the latter a {\it Type II approximate cycle}.  Then we can show:

\begin{theorem} \label{TypeII} Let $F(n)$ be any function so that $q^{-n}F(n) \rightarrow \infty$ as $n \rightarrow \infty$.  Then there exists a $q$-ary Type II approximate $\cI$-cycle of length $F(n)$ when $|\cI|=n$.
\end{theorem}

We need a result of Janson to proceed.  The following appears in \cite{Jan98}.  First, some notation.  Let $I$ be an index set for a set of events $\{B_i\}_{i \in I}$.  Define a graph $\sim$ on $I$ with the following property: Let $J_1$ and $J_2$ be two disjoint subsets of $I$ such that there is no $i_1 \in J_1$ and $i_2 \in J_2$ with $i_1 \sim i_2$.  Now, let $A^1$ be any Boolean function of the events $\{B_i : i \in J_1\}$ and let $A^2$ be any Boolean function of the events $\{B_i : i \in J_2\}$. Then $A^1$ and $A^2$ are independent.

Let $\mu = \sum_{i=1}^m \bP(B_i)$, $\Delta = \sum_{i \sim j} \bP(B_i \wedge B_j)$, and $\delta = \max_i \sum_{j \sim i} \bP(B_j)$.  Then the following holds.

\begin{lemma} With the above notation,
$$
\bP(\wedge_{i=1}^m \bar{B}_i) \leq \exp(-\min\left (\frac{\mu^2}{8\Delta},\frac{\mu}{2},\frac{\mu}{6\delta}\right)).
$$
\end{lemma}

\begin{proof}[Proof of Theorem \ref{TypeII}.] Fix an integer $m$.  We wish to show that (for a suitable choice of $m$) the expected number of $q$-ary $n$-words which do not appear as $\chi(\cI+j)$ for any $j \in \Z_m$ for a random function $\chi : \Z_m \rightarrow q$ is $o(q^n)$.  This quantity is $q^n$ times the probability that a single word -- say, $0^n$ -- does not appear anywhere.  So we need only show that this probability is $o(1)$.  Let $B_j$ be the event that $\chi(\cI+j) = W$.  Then define a graph on the $B_j$ as follows: $B_j \sim B_k$ iff $(\cI+j) \cap (\cI + k) \neq \emptyset$.  Note that $\deg(B_j) \leq n^2$ for all $j$, and $\bP(B_i) = q^{-n}$.  Therefore $\mu = mq^{-n}$, 
$$
\Delta \leq q^n n^2 q^{-n} = n^2
$$
and $\delta \leq n^2 q^{-n}$.  Plugging into the lemma, we find
\begin{align*}
\bP(\wedge_{i=1}^m \bar{B}_i) \leq \exp(-\min\left(\frac{m^2}{8 n^2q^n},\frac{m}{2q^n},\frac{m}{6n^2}\right)).
\end{align*}
If we let $m = F(n)$, then $\min({m^2}/{8 n^2q^n},{m}/{2q^n},{m}/{6n^2}) \rightarrow \infty$, completing the proof.
\end{proof}

Now, we address the problem of constructing Type I approximate $\cI$-cycles.  First, for a set of reals $X = \{x_i\}_{i=1}^n$, define
$$
\mu(X) = \max_{\alpha \in \R} \min_{j \neq k} \| \alpha(x_j-x_k) \|,
$$
where $\|x\|$ denotes the distance to the closest integer.  It is proven in \cite{KRS00} that, for any set $X$ of cardinality $n$, $\mu(X) \geq n^{-2}$.  Then we have the following.

\begin{lemma} \label{fixitup} For any $\cI=\{i_j\}_{j=1}^n \subset \Z$ and collection of $q$-ary $n$-words $W_1,\ldots,W_S$, there exists an integer $p = n^2 S (1 + o(1))$ and a map $\chi : Z_p \rightarrow [q]$ so that every $W_j$ appears as $\chi(\cI+t)$ for some $t \in Z_p$.
\end{lemma}
\begin{proof} Let $p$ be the smallest prime greater than $n^2 (S + 3)$, and choose $\alpha$ achieving the bound $\mu(X) \geq n^{-2}$.  Then there exists a $k \in Z_p$ so that $|\alpha - k/p| \leq p^{-1}$, and it follows that, for any $j \neq k$, $ki_j$ and $ki_k$ are separated by at least $p n^{-2} - 3 \geq S$ integers modulo $p$.  Write $W(k)$ for the $k^\textrm{th}$ symbol of the word $W$.  Then we may define $\chi : \Z_p \rightarrow [q]$ by $\chi(ki_j + t) = W_t(j)$ for $1 \leq t \leq S$ and $1 \leq j \leq n$, since the minimum gap between elements of $k\cI$ is at least $S$.  Define $\chi(s)$ arbitrarily for $s \neq ki_j + t$ for any $j$ and $t$.  Then the map $\chi^\prime : s \rightarrow \chi(k s)$ has the desired property.
\end{proof}

We could use this result immediately to achieve a length $n^2 q^n (1+o(1))$ Type I approximate $\cI$-cycle, but it possible to do better using the random construction above.  There is a trivial lower bound of $q^n$ on the length of any approximate cycle, and we can show an upper bound that is only slightly worse:

\begin{theorem} For any $\cI=\{i_j\}_{j=1}^n \subset \Z$, there exists a $q$-ary Type I approximate $\cI$-cycle of length $p = (8+o(1))q^n \log n $.
\end{theorem}
\begin{proof} The basic idea is to take a random sequence $T_1$ that contains almost all words, then use Lemma \ref{fixitup} to ``tack on'' the remaining ones.  We use the notation of the proof of Theorem \ref{TypeII}.  If we let $F(n) = \ceil{4 q^n \log n}$, then, for sufficiently large $n$, the expected number of words which do not appear as $\chi(\cI+t)$ is at most
$$
q^n \exp(-2 \log n) = \frac{q^n}{n^2}.
$$
Apply the lemma to find a sequence $T_2$ in which each word missed by the random sequence occurs.  Then, we concatenate two copies of $T_1$ with two copies of $T_2$ (two copies are used to avoid disturbing sequences which ``wrap around''), and the result is a Type I approximate $\cI$-cycle of the stated length.
\end{proof}

\section{Decompositions into Directed Cycles} \label{sec:trails}

A number of recent papers (most notably \cite{AGSV03},\cite{B01}, and \cite{B03}) have addressed (and solved) the problem of decomposing a complete (possibly directed) graph into a set of cycles of prescribed length.  Generally speaking, so long as the cycle lengths add up to the number of edges, there are very few impediments to the existence of such decompositions -- although demonstrating this is far from simple.  None of this work has dealt with graphs containing loops, however; in this section, we address this situation, in the case when all the cycles have the same length.

Let $\KD_n$ denote the complete directed graph with loops on $n$ vertices, i.e., the vertex set is $[n]$ and the edge set is $[n] \times [n]$, and let $\KL_n$ denote the complete directed graph without loops.  We wish to know, for which $n$ and $d$ is it possible to decompose the edge set of $\KD_n$ completely into cycles of length $d$?  Clearly $d$ must divide $n^2$.  However, our main result states that the only other obstruction is that $d > 2$.  (Indeed, it is easy to see in these cases that such a decomposition is not possible.)

First, we state a result from \cite{B03}.

\begin{theorem} \label{prevresult} If $\sum^t_{i=1} m_i = n(n - 1)$ and $m_i \geq 2$ for $i = 1, \ldots , t$, then $\KL_n$ can be decomposed as the edge-disjoint union of cycles of lengths $m_1, \ldots , m_t$, except in the case when $n=6$ and all $m_i = 3$.
\end{theorem}

This will imply the following:

\begin{prop} If $d | n^2$, $d=6$ or $d \geq 8$, then $\KD_n$ may be decomposed into cycles of length $d$.
\end{prop}

\begin{proof} We offer a procedure for ``packing'' length $d$ cycles first into $G$, then the $G_j$'s, then $H$.  Split the vertex set of $\KD_n$ into two pieces: $U = \{a,b\}$ and $V = \{1,\ldots,n-2\}$.  We may then decompose the edge set of $\KD$ into the following pieces:
\begin{enumerate}
\item one $\KD_2$ on $U$,
\item one $\KL_{n-2}$ on $V$,
\item and $n-2$ graphs each of which has vertex set $\{a,b,j\}$ for some $j \in V$, with a loop at $j$ and edges $(j,i)$ and $(i,j)$ for $i = a,b$.
\end{enumerate}
Such a decomposition is possible because $d \geq 6$ implies $n > 2$.

Call the first graph $H$, the second $G$, and the third $G_j$.  Denote by $\{u,v\}$ the pair $\{(u,v),(v,u)\}$.  Suppose $\binom{n-2}{2} = r \! \mod d$ with $1 \leq r \leq d$.  Supposing $r > 1$, by Theorem \ref{prevresult}, we may decompose $G$ into $K=\floor{(\binom{n-2}{2}-1)/d}$ length $d$ cycles and one $r$-trail, which we call $T$.  We may assume, without loss of generality, that $T$ meets vertex $3$.  Let $m = \floor{(d-r)/5}$ and $d^\prime = d-r-5m$.  We construct a graph $X$ as follows.  Take $X$ to be the union of $T$, $G_j$ for $j=1, \ldots m-1$, and one of the following, according to the value of $d^\prime$.  Note that $d^\prime < 4$ implies $m \leq n-3$, since otherwise the number of remaining edges (i.e, ones unaccounted for thus far by $X$ or any of the cycles in $\KL_{n-2}$) would not be divisible by $d$.  Similarly, $d^\prime = 4$ implies either $m < n-3$ or we may add in all the remaining edges of the graph and be done.
\begin{itemize}
\item $d^\prime = 0$ : Add $G_m$ to $X$.
\item $d^\prime = 1$ : If $m=0$, add the edge $(1,1)$ to $X$.  If $m > 0$, add the edges $\{a,m\}$, $(m,m)$, $\{a,m+1\}$, and $(m+1,m+1)$.
\item $d^\prime = 2$ : Add $\{a,m+1\}$.
\item $d^\prime = 3$ : Add $\{a,m+1\}$ and $(m+1,m+1)$.
\item $d^\prime = 4$ : Add $\{a,m+1\}$ and $\{a,m+2\}$.
\end{itemize}
Note that in all cases, $X$ is connected and has equal in- and out-degree at every vertex; therefore, $X$ is Eulerian and may be written as a cycle of length $d$.

Now, if not all edges of $\KD_n$ have been accounted for, yet, we wish to add another length $d$ cycle.  At the previous step, there are a few possibilities for the remaining set of edges not assigned to cycles.  The set contains $H$, some of the $G_j$, as well as either:
\begin{enumerate}
\item Case I: $\{a,m\},\{b,m\}$,
\item Case II: $\{b,m\},\{b,m+1\}$,
\item Case III: $(m,m),\{b,m\}$,
\item Case IV: $\{b,m\}$, or
\item Case V: $(m,m),\{b,m\},(m+1,m+1),\{b,m+1\}$.
\end{enumerate}
In each of the cases, let $Y$ be the set of edges listed above.  We now let $m = \floor{(d - |Y|)/5}$ and $d^\prime = d - 5m - |Y|$.  We may proceed exactly as above with the construction of $X$, unless $m=0$.  If $m=0$, then either $d^\prime \geq 2$ or $d^\prime = 0$.  In the latter case, $d = 6$ and we are in Case V, so we may simply take $X = Y$.  In the former case, we append $Y$, as well as $d^\prime > 0$ edges drawn from the first full ``unused'' $G_k$ as follows.  (Note that, if there is no unused $G_k$, we may append the remaining edges in the graph and be done.)
\begin{itemize}
\item $d^\prime = 2$ : Add $\{b,k\}$.
\item $d^\prime = 3$ : Add $\{b,k\}$ and $(k,k)$.
\item $d^\prime = 4$ : If $k < n-2$, add $\{b,k\}$ and $\{b,k+1\}$. It is not possible for $k=n-2$, since there would be too few edges left for $d$ to divide them evenly.
\end{itemize}

Again, the resulting graph $X$ is Eulerian, and can therefore be written as a length $d$ cycle.  It is clear that we may repeat the previous step until a full decomposition is achieved, though we may need to switch the roles of $a$ and $b$.

It remains to deal with the case of $r=1$ at the beginning of the proof.  Instead of decomposing into all length $d$ trails except for one $r$-trail, we instead create one length $d-1$ trail and one length $2$ trail.  Then attach $(a,a)$ to the $d-1$ trail (changing labels if necessary to maintain connectivity), and proceed as above with $r=2$, using the remaining edges of $H$ when necessary.  All of the details of the preceding argument work with this slight modification.
\end{proof}

\begin{prop} If $n$ is even, then $\KD_n$ may be decomposed into cycles of length $4$.
\end{prop}
\begin{proof} There are two cases: $n = 2 \! \mod 4$ and $n = 0 \! \mod 4$.  Suppose the former.  Then the following types of edge-sets partition $E(\KD_n)$ into cycles of length $4$:
\begin{enumerate}
\item $\{(j,j), (j+n/2,j+n/2), (j,j+n/2), (j+n/2,j)\}$ for each $j$ with $1 \leq j \leq n/2$,
\item $\{(j,j+2k-1),(j+2k-1,j),(j,j+2k),(j+2k,j)\}$, for each $j$ and each $k$ so that $1 \leq k \leq (n-2)/4$.
\end{enumerate}
If $4 | n$, we use the following types of sets instead:
\begin{enumerate}
\item $\{(j,j), (j+n/2,j+n/2), (j,j+n/2), (j+n/2,j)\}$ for each $j$ with $1 \leq j \leq n/2$,
\item $\{(j,j+2k),(j+2k,j),(j,j+2k+1),(j+2k+1,j)\}$, for each $j$ and each $k$ so that $1 \leq k \leq n/4-1$.
\item $\{(2j,2j-1),(2j-1,2j),(2j,2j+1),(2j+1,2j)\}$, for each $j$ with $1 \leq j \leq n/2$.
\end{enumerate}
\end{proof}

\begin{prop} If $d | n$, $d = 3,5,$ or $7$, then $\KD_n$ may be decomposed into cycles of length $d$.
\end{prop}
\begin{proof} We imitate the case of $d \geq 8$ here, only the situation is simpler.  Suppose $d=3$.  We may decompose $\KD_n$ into three pieces: the loop $(a,a)$, a $\KL_{n-1}$, and $n-1$ length $3$ cycles $G_1,\ldots,G_{n-1}$.  Applying Theorem \ref{prevresult}, we may break the second of these into cycles of length $3$, except for one of length $2$, which we call $T$.  (Since $3 | n^2$, $(n-1)(n-2) = 2 \! \mod 3$.)  We may assume that $T$ meets $a$; therefore, appending $(a,a)$ to $T$ and taking each $G_j$ as its own cycle provides a decomposition.

Now suppose $d=5$.  We may decompose $\KD_n$ into three pieces again: a $\KD_2$ on $\{a,b\}$, a $\KL_{n-2}$, and $n-2$ length $5$ cycles $G_1,\ldots,G_{n-2}$.  Applying Theorem \ref{prevresult}, we may break the $\KL_{n-2}$ into cycles of length $5$, except for one of length $4$ and one of length $2$, which we call $T_1$ and $T_2$, respectively.  (This time, $(n-2)(n-3) = 1 \! \mod 5$.)  We may assume that $T_1$ meets $a$ and $T_2$ meets $b$.  Append $(a,a)$ to $T_1$, $\KD_2 \setminus (a,a)$ to $T_2$, and take each $G_j$ as its own element of the decomposition.

Finally, let $d=7$.  We have the following decomposition of $\KD_n$: a $\KD_3$ on $\{a,b,c\}$, a $\KL_{n-3}$, and $n-3$ length $7$ cycles $G_1,\ldots,G_{n-3}$.  Applying Theorem \ref{prevresult}, we may break the $\KL_{n-3}$ into cycles of length $7$, except for one of length $5$, which we call $T^\prime$.  (Now, $(n-3)(n-4) = 5 \! \mod 7$.)  We may assume that $T^\prime$ meets $a$.  Append $\{a,b\}$ to $T^\prime$, include the $7$-cycle $\KD_3 \setminus \{a,b\}$, and take each $G_j$ as its own element of the decomposition.
\end{proof}

All of these results together imply the following.

\begin{theorem} \label{trails} If $d | n^2$ and $d \geq 3$, then $\KD_n$ may be decomposed into cycles of length $d$.
\end{theorem}

\section{Conclusion}

We wish to know, first and foremost, what distinguishes valid sets from invalid ones.  The authors' attempts to find a simple way to separate these cases was met with frustration.  A simpler problem is the case of index sets which are arithmetic progressions: the question of decomposing de Bruijn graphs into cycles is a natural one.  We would also like to see corresponding lower bounds or improvements on the bounds of Section \ref{sec:approx}, particularly in the case of Type I approximate cycles.

Finally, the problem solved by Theorem \ref{trails} has a natural generalization along the lines of other, similar work.  Suppose that $\{m_j\}$ is such that $\sum_j m_j = n^2$.  When is it possible to decompose $\KD_n$ into cycles of lengths $\{m_j\}$?  Clearly, there can be at most $n$ 1's among the $m_j$, and there must be at least $n$ indices $j$ for which $m_j \neq 2$.  Are there additional restrictions?

\section{Acknowledgements}

Thank you to Paul Balister for the helpful suggestion concerning Theorem \ref{trails}.

\end{document}